\documentclass[11pt]{amsart}
 \usepackage{amssymb,amsfonts,amscd,amsbsy,pst-plot,pst-node}
\usepackage[mathscr]{eucal}
\usepackage{url}

\theoremstyle{plain}
\newtheorem{thm}{Theorem}[section]
\newtheorem{lem}[thm]{Lemma}
\newtheorem{prop}[thm]{Proposition}

\theoremstyle{definition}

\numberwithin{equation}{section}

\begin{document}
\title[$M_2$-rank differences for partitions without repeated odd parts]{$M_2$-rank differences for partitions without repeated odd parts}
\author{Jeremy Lovejoy and Robert Osburn}

\address{CNRS, LIAFA, Universit{\'e} Denis Diderot, 2, Place Jussieu, Case 7014, F-75251 Paris Cedex 05, FRANCE}

\address{School of Mathematical Sciences, University College Dublin, Belfield, Dublin 4, Ireland}

\address{IH{\'E}S, Le Bois-Marie, 35, route de Chartres, F-91440 Bures-sur-Yvette, FRANCE}

\email{lovejoy@liafa.jussieu.fr}

\email{robert.osburn@ucd.ie, osburn@ihes.fr}

\thanks{The first author was partially supported by an ACI ``Jeunes Chercheurs et Jeunes Chercheuses".}


\date{May 31, 2007}

\begin{abstract}
We prove formulas for the generating functions for $M_2$-rank
differences for partitions without repeated odd parts. These
formulas are in terms of modular forms and generalized Lambert
series.
\end{abstract}

\maketitle
\section{Introduction}
A partition of a non-negative integer $n$ is a non-increasing
sequence whose sum is $n$.  One of the most useful ways to
represent a partition is with the Ferrers diagram.  For example,
the partition $(10,6,6,3,1)$ is represented by the diagram \\
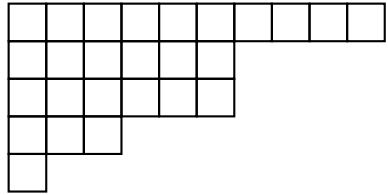
\begin{center}
    \psset{unit=0.5}
    \begin{pspicture}(0,0)(5,-5)
    \psframe[dimen=middle](0,0)(1,-1)
    \psframe[dimen=middle](1,0)(2,-1)
    \psframe[dimen=middle](2,0)(3,-1)
    \psframe[dimen=middle](3,0)(4,-1)
    \psframe[dimen=middle](4,0)(5,-1)
    \psframe[dimen=middle](5,0)(6,-1)
    \psframe[dimen=middle](6,0)(7,-1)
    \psframe[dimen=middle](7,0)(8,-1)
    \psframe[dimen=middle](8,0)(9,-1)
    \psframe[dimen=middle](9,0)(10,-1)
    \psframe[dimen=middle](0,-1)(1,-2)
    \psframe[dimen=middle](1,-1)(2,-2)
    \psframe[dimen=middle](2,-1)(3,-2)
    \psframe[dimen=middle](3,-1)(4,-2)
    \psframe[dimen=middle](4,-1)(5,-2)
    \psframe[dimen=middle](5,-1)(6,-2)
    \psframe[dimen=middle](0,-2)(1,-3)
    \psframe[dimen=middle](1,-2)(2,-3)
    \psframe[dimen=middle](2,-2)(3,-3)
    \psframe[dimen=middle](3,-2)(4,-3)
    \psframe[dimen=middle](4,-2)(5,-3)
    \psframe[dimen=middle](5,-2)(6,-3)
    \psframe[dimen=middle](0,-3)(1,-4)
    \psframe[dimen=middle](1,-3)(2,-4)
    \psframe[dimen=middle](2,-3)(3,-4)
    \psframe[dimen=middle](0,-4)(1,-5)
\end{pspicture}
\end{center}
MacMahon \cite{Ma1} generalized the Ferrers diagram to an
$M$-modular diagram of a partition.  A special case of his
construction, the $2$-modular diagram is a Ferrers diagram where
all of the boxes are filled with $2$'s except possibly the last
box of a row, which may be filled with a $1$, with the condition
that no $2$ occurs directly below a $1$. As an illustration, the
partition $(10,10,8,7,7,4,2,2,1)$ has $2$-modular diagram \\

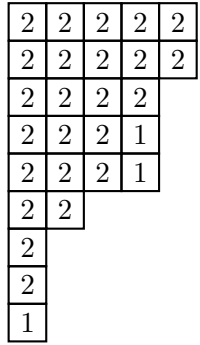
\begin{center}
    \psset{unit=0.5}
    \begin{pspicture}(0,0)(5,-9)
    \psframe[dimen=middle](0,0)(1,-1)
    \psframe[dimen=middle](1,0)(2,-1)
    \psframe[dimen=middle](2,0)(3,-1)
    \psframe[dimen=middle](3,0)(4,-1)
    \psframe[dimen=middle](4,0)(5,-1)
    \rput(0.5,-0.5){2}
    \rput(1.5,-0.5){2}
    \rput(2.5,-0.5){2}
    \rput(3.5,-0.5){2}
    \rput(4.5,-0.5){2}
    \psframe[dimen=middle](0,-1)(1,-2)
    \psframe[dimen=middle](1,-1)(2,-2)
    \psframe[dimen=middle](2,-1)(3,-2)
    \psframe[dimen=middle](3,-1)(4,-2)
    \psframe[dimen=middle](4,-1)(5,-2)
    \rput(0.5,-1.5){2}
    \rput(1.5,-1.5){2}
    \rput(2.5,-1.5){2}
    \rput(3.5,-1.5){2}
    \rput(4.5,-1.5){2}
    \psframe[dimen=middle](0,-2)(1,-3)
    \psframe[dimen=middle](1,-2)(2,-3)
    \psframe[dimen=middle](2,-2)(3,-3)
    \psframe[dimen=middle](3,-2)(4,-3)
    \rput(0.5,-2.5){2}
    \rput(1.5,-2.5){2}
    \rput(2.5,-2.5){2}
    \rput(3.5,-2.5){2}
    \psframe[dimen=middle](0,-3)(1,-4)
    \psframe[dimen=middle](1,-3)(2,-4)
    \psframe[dimen=middle](2,-3)(3,-4)
    \psframe[dimen=middle](3,-3)(4,-4)
    \rput(0.5,-3.5){2}
    \rput(1.5,-3.5){2}
    \rput(2.5,-3.5){2}
    \rput(3.5,-3.5){1}
    \psframe[dimen=middle](0,-4)(1,-5)
    \psframe[dimen=middle](1,-4)(2,-5)
    \psframe[dimen=middle](2,-4)(3,-5)
    \psframe[dimen=middle](3,-4)(4,-5)
    \rput(0.5,-4.5){2}
    \rput(1.5,-4.5){2}
    \rput(2.5,-4.5){2}
    \rput(3.5,-4.5){1}
    \psframe[dimen=middle](0,-5)(1,-6)
    \psframe[dimen=middle](1,-5)(2,-6)
    \rput(0.5,-5.5){2}
    \rput(1.5,-5.5){2}
    \psframe[dimen=middle](0,-6)(1,-7)
    \rput(0.5,-6.5){2}
    \psframe[dimen=middle](0,-7)(1,-8)
    \rput(0.5,-7.5){2}
    \psframe[dimen=middle](0,-8)(1,-9)
    \rput(0.5,-8.5){1}
\end{pspicture}
\end{center}

If we add the condition that a $1$ may only occur in the last
entry of a column, then these $2$-modular diagrams correspond to
partitions whose odd parts may not be repeated.  Partitions
without repeated odd parts and their $2$-modular diagrams have
long played a role in combinatorial studies of $q$-series
identities (see \cite{An1,bg,Be1,Co-Ma1,Pa1,Ze1}, for example).
Most recently, Berkovich and Garvan \cite{bg} introduced what they
called the $M_2$-rank of such partitions. The $M_2$-rank of a
partition $\lambda$ without repeated odd parts is defined to be
the number of columns minus the number of rows of its $2$-modular
diagram, or equivalently,

$$
\text{$M_{2}$-rank} \hspace{.025in}
(\lambda)=\Bigg\lceil{\frac{l(\lambda)}{2}}\Bigg\rceil -
\nu(\lambda),
$$

\noindent where $l(\lambda)$ is the largest part of $\lambda$ and
$\nu(\lambda)$ is the number of parts of $\lambda$.

The two-variable generating function for the $M_2$-rank has a
particularly nice form.  Namely, it may be deduced from
\cite[Theorem 1.2]{Lo1} that if $N_2(m,n)$ denotes the number of
partitions of $n$ without repeated odd parts whose $M_2$-rank is
$m$, then

\begin{equation} \label{gf1}
\sum_{n \geq 0 \atop m \in \mathbb{Z}} N_2(m,n)z^mq^n
=\sum_{n=0}^{\infty} q^{n^2} \frac{(-q; q^2)_{n}}{(zq^2,
q^2/z)_{n}}.
\end{equation}
Here we have introduced the standard $q$-series notation
\cite{Ga-Ra1},
$$
(a_1,a_2,\dots,a_j;q)_n = \prod_{k=0}^{n-1}(1-a_1q^k)(1-a_2q^k)
\cdots (1-a_jq^k),
$$
following the custom of dropping the ``$;q$" unless the base is
something other than $q$.

It has recently been discovered that the generating function for
the $M_2$-rank in \eqref{gf1} has a nice number-theoretic
structure.  Bringmann, Ono, and Rhoades \cite{Br-On-Rh1} have
shown that one obtains the holomorphic part of a weak Maass form
when $z$ is replaced by certain roots of unity.  When $z=i$, this
is the eighth order mock theta function $U_0(q)$ of Gordon and
McIntosh \cite{Go-Mc1} (c.f. \cite[p.29]{Lost}). There are many
nice consequences of this number-theoretic structure, including
the fact that the generating function for $N_2(s,\ell,n) -
N_2(t,\ell,n)$ will often be a classical modular form when $n$ is
restricted to arithmetic progressions. Here $N_{2}(s, \ell, n)$
denotes the number of partitions of $n$ without repeated odd parts
whose $M_{2}$-rank is congruent to $s$ modulo $\ell$. In this
paper we obtain formulas for all of the generating functions
$N_2(s,\ell,\ell n+d)-N_2(t,\ell,\ell n+d)$, when $\ell = 3$ or
$5$, in terms of modular forms and generalized Lambert series. We
shall indeed see that many of these functions are simply modular
forms.

Using the notation
\begin{equation} \label{rst}
R_{st}(d) = \sum_{n \geq 0} \left({N}_2(s,\ell,\ell n + d) -
{N}_2(t,\ell,\ell n + d)\right ) q^n,
\end{equation}
where the prime $\ell$ will always be clear, the main results are
summarized in Theorems \ref{main3} and \ref{main5} below.

\begin{thm} \label{main3} For $\ell=3$, we have

\begin{equation} \label{r_01(0)}
R_{01}(0) =-1 -3q^3 \frac{(-q^3; q^{6})_{\infty}}{(q^{6};
q^{6})_{\infty}} \sum_{n \in \mathbb{Z}}
\frac{(-1)^nq^{6n^2+9n}}{1-q^{6n+4}}
 + \frac{(q^{6}; q^{6})_{\infty}^4 (-q^3; q^3)_{\infty}^{4} (q; q^2)_{\infty}}
 {(q^{4}; q^{4})_{\infty} (q^2, q^{10}, q^{12}; q^{12})_{\infty}^2},
\end{equation}

\begin{equation} \label{r_01(1)}
R_{01}(1)=\frac{(-q^3, q^6; q^6)_{\infty}}{(q^2, q^4; q^6)_{\infty}},
\end{equation}

\begin{equation} \label{r_01(2)}
R_{01}(2)=\frac{(q^3; q^3)_{\infty} (-q^6; q^6)_{\infty}}{(q, q^5; q^6)_{\infty} (q^4, q^8; q^{12})_{\infty}}.
\end{equation}

\end{thm}

\begin{thm} \label{main5} For $\ell=5$, we have

\begin{equation} \label{r_12(0)}
\begin{aligned}
R_{12}(0) & =-1 -q^2 \frac{(-q^{5}; q^{10})_{\infty}}{(q^{10}; q^{10})_{\infty}}
\sum_{n \in \mathbb{Z}}\frac{(-1)^nq^{10n^2+15n}}{1-q^{10n+2}}  \\
& + \frac{(q, q^{9}; q^{10})_{\infty}^2 (q^{6}, q^{8}, q^{12}, q^{14}; q^{20})_{\infty}
 (q^{10}; q^{20})_{\infty}^3 (q^{20}; q^{20})_{\infty}^2}{(q; q)_{\infty}},
\end{aligned}
\end{equation}

\begin{equation} \label{r_12(1)}
R_{12}(1)=0,
\end{equation}

\begin{equation} \label{r_12(2)}
R_{12}(2)=\frac{q(q^2, q^{18}; q^{20})_{\infty} (q^5; q^5)_{\infty} (-q^{10}; q^{10})_{\infty}}{(q, q^4; q^5)_{\infty}},
\end{equation}

\begin{equation} \label{r_12(3)}
R_{12}(3)=\frac{(-q^5, q^{10}; q^{10})_{\infty}}{(q^4, q^6; q^{10})_{\infty}},
\end{equation}

\begin{equation} \label{r_12(4)}
R_{12}(4)=2q^3 \frac{(-q^{5}; q^{10})_{\infty}}{(q^{10};
q^{10})_{\infty}} \sum_{n \in \mathbb{Z}}
\frac{(-1)^nq^{10n^2+15n}}{1-q^{10n+4}} + \frac{(q^{3}, q^{7},
q^{10}; q^{10})_{\infty}^2} {(q; q^{2})_{\infty} (q^{6}, q^{8},
q^{12}, q^{14}, q^{20}; q^{20})_{\infty}},
\end{equation}

\begin{equation} \label{r_02(0)}
R_{02}(0)=1 + 2q^2  \frac{(-q^{5}; q^{10})_{\infty}}{(q^{10};
q^{10})_{\infty}} \sum_{n \in
\mathbb{Z}}\frac{(-1)^nq^{10n^2+15n}}{1-q^{10n+2}} - \frac{(q,
q^{9}; q^{10})_{\infty}^2 (q^{10}; q^{10})_{\infty}^3 (q^{6},
q^{8}, q^{12}, q^{14}; q^{20})_{\infty}} {(q; q)_{\infty} (q^{20};
q^{20})_{\infty}},
\end{equation}

\begin{equation} \label{r_02(1)}
R_{02}(1)=\frac{(-q^5, q^{10}; q^{10})_{\infty}}{(q^2, q^8; q^{10})_{\infty}},
\end{equation}

\begin{equation} \label{r_02(2)}
R_{02}(2)=\frac{(q^5; q^5)_{\infty} (-q^{10}; q^{10})_{\infty} (q^6, q^{14}; q^{20})_{\infty}}{(q^2, q^3; q^5)_{\infty}},
\end{equation}

\begin{equation} \label{r_02(3)}
R_{02}(3)=0,
\end{equation}

\begin{equation} \label{r_02(4)}
R_{02}(4)=q^3 \frac{(-q^{5}; q^{10})_{\infty}}{(q^{10};
q^{10})_{\infty}} \sum_{n \in
\mathbb{Z}}\frac{(-1)^nq^{10n^2+15n}}{1-q^{10n+4}} + \frac{(q^{3},
q^{7}, q^{10}; q^{10})_{\infty}^2}{(q; q^{2})_{\infty} (q^{6},
q^{8}, q^{12}, q^{14}, q^{20}; q^{20})_{\infty}}.
\end{equation}

\end{thm}

To prove Theorems \ref{main3} and \ref{main5}, we shall roughly
follow the method developed by Atkin and Swinnerton-Dyer
\cite{asd} in their study of Dyson's rank for partitions. This method may be generally described as
regarding groups of identities as equalities between polynomials
of degree $\ell-1$ in $q$ whose coefficients are power series in
$q^{\ell}$. Specifically, we first consider the expression

\begin{equation} \label{term}
\sum_{n=0}^{\infty} \Bigl\{ N_{2}(s,\ell,n) -
N_{2}(t,\ell,n) \Bigr\} q^{n}
\frac{(q^2; q^2)_{\infty}}{(-q; q^2)_{\infty}}.
\end{equation}

\noindent By (\ref{gf2}), (\ref{s}), and (\ref{final}), we write
(\ref{term}) as a polynomial in $q$ whose coefficients are power
series in $q^{\ell}$.  We then alternatively express (\ref{term})
in the same manner using Theorems \ref{main3} and \ref{main5} and Lemma \ref{lem6}.
Finally, we use various $q$-series identities to show that these
two resulting polynomials are the same for each pair of values of
$s$ and $t$.

The paper is organized as follows.  In Section 2 we collect some
basic definitions, notations and generating functions.  In Section
3 we record a number of equalities between an infinite product and
a sum of infinite products.  These are ultimately required for the
simplification of identities that end up being more complex than
we would like, principally because there are only two $0$'s in
Theorems \ref{main3} and \ref{main5}. In Section 4 we prove two key $q$-series
identities relating generalized Lambert series to infinite
products, and in Section 5 we give the proofs of Theorems
\ref{main3} and \ref{main5}.

\section{Preliminaries}
 We begin by introducing some notation and definitions, essentially
following \cite{asd}. With $y=q^{\ell}$, let

$$
r_{s}(d):=\sum_{n=0}^{\infty} N_{2}(s, \ell,  \ell n+d) y^n
$$

\noindent and

$$
r_{st}(d):=r_{s}(d) - r_{t}(d).
$$

\noindent Thus we have

$$
\sum_{n=0}^{\infty} N_{2}(s, \ell ,n) q^n = \sum_{d=0}^{\ell-1}
r_{s}(d) q^d.
$$

To abbreviate the sums appearing in Theorems \ref{main3} and \ref{main5}, we define

$$
\Sigma(z,\zeta, q):= \sum_{n \in \mathbb{Z}} \frac{(-1)^n
\zeta^{4n} q^{2n^2 + 3n}}{1-z^{2}q^{2n}}.
$$

\noindent Henceforth we assume that $a$ is not a multiple of $q$.  We
write
$$
\Sigma(a,b):=\Sigma(y^a, y^b, y^{\ell})=\sum_{n \in \mathbb{Z}}
\frac{(-1)^n y^{4bn + {\ell}n(2n+3)}}{1-y^{2{\ell}n+2a}}
$$

\noindent and

$$
\Sigma(0,b) :=  \sideset{}{'} \sum_{n \in \mathbb{Z}} \frac{(-1)^n
y^{4bn + \ell n(2n+3)}}{1-y^{2\ell n}},
$$
where the prime means that the term corresponding to $n=0$ is
omitted.

To abbreviate the products occurring in Theorems \ref{main3} and \ref{main5}, we
define

$$
P(z,q):=\prod_{r=1}^{\infty} (1-zq^{r-1})(1-z^{-1}q^r)
$$

%

\noindent and

$$
P(0):=\prod_{r=1}^{\infty} (1-y^{2\ell r}).
$$


\noindent We have the relations

\begin{equation} \label{p1}
P(z^{-1}q, q)=P(z,q)
\end{equation}

\noindent and

\begin{equation} \label{p2}
P(zq, q)=-z^{-1}P(z,q).
\end{equation}





Now, for any integer
$m$ we have \cite{Lo1}
\begin{equation} \label{gf1.5}
\sum_{n \geq 0} N_2(m,n)q^n =
\frac{(-q;q^2)_{\infty}}{(q^2;q^2)_{\infty}}\sum_{n \geq 1}
(-1)^{n+1}q^{2n^2-n + 2|m|n}(1-q^{2n}).
\end{equation}
It is then a simple matter to deduce that
\begin{equation} \label{gf2}
\sum_{n=0}^{\infty} N_{2}(s, m, n) q^n = \frac{(-q;
q^2)_{\infty}}{(q^2; q^2)_{\infty}} \sideset{}{'} \sum_{n \in
\mathbb{Z}} \frac{(-1)^n q^{2n^2 + n}(q^{2sn} + q^{2(m-s)n})}{1 -
q^{2mn}}.
\end{equation}
Unfortunately, it does not appear that one can go
directly from differences of \eqref{gf2} to the formulas in
Theorems \ref{main3} and \ref{main5}. Hence it will be beneficial to consider sums of the form

\begin{equation} \label{s}
S_{2}(b):= \sideset{}{'}\sum_{n \in \mathbb{Z}}
\frac{(-1)^n q^{2n^2+ bn}}{1-q^{2 \ell n}}.
\end{equation}

\noindent We will require the relation

\begin{equation} \label{rels}
S_{2}(b)=-S_{2}(2 \ell -b),
\end{equation}

\noindent which follows from the substitution $n \to -n$ in (\ref{s}). We
shall also require the fact that

\begin{equation} \label{bodd}
S_{2}(b)-S_{2}(2\ell +b) = \sum_{n \in \mathbb{Z}} (-1)^n q^{2n^2 + bn} - 1 =
(q^{2+b}, q^{2-b}, q^4; q^4)_{\infty} -1
\end{equation}

\noindent if $b$ is odd. This follows by applying the case $z=-q^b$ and replacing $q$
with $q^2$ in Jacobi's triple product identity

\begin{equation} \label{jtp}
\sum_{n \in \mathbb{Z}} z^nq^{n^2} = (-zq,-q/z,q^2;q^2)_{\infty}.
\end{equation}

\section{Infinite Product Identities}
In this section we record some identities involving infinite
products.  These will be needed later on for simplification and
verification of certain identities. First, we have a result which
is the analogue of Lemma 6 in \cite{asd} and Lemma 3.1 in
\cite{loveoz}. The proof, which just amounts to an application of
\eqref{jtp}, is similar to that of Lemma 3.1 in \cite{loveoz} and
thus is omitted.

\begin{lem} \label{lem6} We have

\begin{equation} \label{lem6eq1}
\frac{(q^2; q^2)_{\infty}}{(-q; q^2)_{\infty}} = (q^3, -q^6, -q^9, -q^{12}, q^{15}, q^{18}; q^{18})_{\infty} -
q(q^9, q^{27}, q^{36}; q^{36})_{\infty}
\end{equation}
and
\begin{equation} \label{lem6eq2}
\begin{aligned}
\frac{(q^2; q^2)_{\infty}}{(-q; q^2)_{\infty}} &= (-q^{10}, q^{15}, -q^{25}, q^{35}, -q^{40}, q^{50}; q^{50})_{\infty} \\
& - q(q^5, -q^{20}, -q^{25}, -q^{30}, q^{45}, q^{50}; q^{50})_{\infty} - q^{3}(q^{25}, q^{75}, q^{100}; q^{100})_{\infty}. \\
\end{aligned}
\end{equation}

\end{lem}

Next, we quote a result of Hickerson \cite[Theorem 1.1]{Hi1} along
with some of its corollaries.

\begin{lem} \label{Hick1}
$$
P(x,q)P(z,q)(q)_{\infty}^2 =
P(-xz,q^2)P(-qz/x,q^2)(q^2;q^2)_{\infty}^2 -
xP(-xzq,q^2)P(-z/x,q^2)(q^2;q^2)_{\infty}^2.
$$
\end{lem}

The first corollary was recorded by Hickerson \cite[Theorem
1.2]{Hi1}.  It follows by applying Lemma \ref{Hick1} twice, once
with $x$ replaced by $-x$ and once with $z$ replaced by $-z$, and
then subtracting.
\begin{lem} \label{Hick2}
$$
P(-x,q)P(z,q)(q)_{\infty}^2 - P(x,q)P(-z,q)(q)_{\infty}^2 =
2xP(z/x,q^2)P(xzq,q^2)(q^2;q^2)_{\infty}^2.
$$
\end{lem}
The second corollary follows just as the first, except we add
instead of subtract in the final step.
\begin{lem} \label{Hick2.5}
$$
P(-x,q)P(z,q)(q)_{\infty}^2 + P(x,q)P(-z,q)(q)_{\infty}^2 =
2P(xz,q^2)P(qz/x,q^2)(q^2;q^2)_{\infty}^2.
$$
\end{lem}


\section{Two Lemmas}

Theorems \ref{main3} and \ref{main5} will follow from identities which relate the sums $\Sigma(a,b)$ to products $P(z,w)$. The key steps are the two Lemmas below. These results are similar in nature to Lemmas 7 and 8 in \cite{asd} and Lemmas 4.1 and 4.2 in \cite{loveoz}.

\begin{lem} \label{chan} We have

\begin{equation} \label{chaneq}
\begin{aligned}
& \sum_{n=-\infty}^{\infty} (-1)^n q^{2n^2 + 3n} \Bigl[ \frac{\zeta^{-4n}}{1-z^{2}{\zeta^{-2}}q^{2n}} +
\frac{\zeta^{4n+6}}{1-z^{2}{\zeta^{2}}q^{2n}} \Bigr] \\
& = \frac{\zeta^{2} (-q, -q, \zeta^{4}, q^{2} \zeta^{-4}; q^2)_{\infty}}{(\zeta^{2}, q^{2} \zeta^{-2}, -q\zeta^{2}, -q\zeta^{-2}; q^{2})_{\infty}} \sum_{n=-\infty}^{\infty} (-1)^n \frac{q^{2n^2 + 3n}}{1-z^{2}q^{2n}} \\
&+ \frac{(-z^{2}q, -qz^{-2}, \zeta^{4}, q^{2} \zeta^{-4}, \zeta^{2}, q^{2} \zeta^{-2}; q^2)_{\infty} (q^2; q^2)_{\infty}^{2}}{(z^2 \zeta^{-2}, q^2 \zeta^2 z^{-2}, z^2 \zeta^2, q^2 z^{-2} \zeta^{-2}, z^2, q^2 z^{-2}, -q\zeta^{2}, -q\zeta^{-2}; q^2)_{\infty}}.
\end{aligned}
\end{equation}

\end{lem}

\begin{proof}
This is just the case $r=1$, $s=3$, $q=q^2$, $a_1=-z^2q$,
$b_1=z^2/\zeta^2$, $b_2 = z^2\zeta^2$, and $b_3=z^2$ of
\cite[Theorem 2.1]{Ch1},
\begin{equation}
\begin{aligned}
& \frac{P(a_1,q)\cdots P(a_r,q) (q)_{\infty}^2}{P(b_1,q) \cdots
P(b_s,q)} \\ & = \frac{P(a_1/b_1,q)\cdots
P(a_r/b_1,q)}{P(b_2/b_1,q)\cdots P(b_s/b_1,q)} \sum_{n \in
\mathbb{Z}}\frac{(-1)^{(s-r)n}q^{(s-r)n(n+1)/2}}{1-b_1q^n}
\left(\frac{a_1\cdots a_r b_1^{s-r-1}}{b_2\cdots b_s}\right)^n \\
&+ \text{idem}(b_1;b_2,\dots,b_s).
\end{aligned}
\end{equation}
Here we use the usual notation
$$
\begin{aligned}
F(b_1,b_2,\dots ,b_m) &+  \text{idem}(b_1;b_2,\dots ,b_m) \\ &:=
F(b_1,b_2, \dots ,b_m) + F(b_2,b_1,b_3, \dots ,b_m) + \cdots +
F(b_m,b_2,\dots,b_{m-1},b_1).
\end{aligned}
$$
\end{proof}

We now specialize Lemma \ref{chan} to the case $\zeta=y^a$, $z=y^b$, and $q=y^{\ell}$:

\begin{equation} \label{lem1}
\begin{aligned}
y^{6a} \Sigma(b+a, a) + \Sigma(b-a, -a) -&  y^{2a} \frac{P(-y^{\ell}, y^{2\ell})P(y^{4a}, y^{2\ell})}
{P(y^{2a}, y^{2\ell})P(-y^{2a+\ell}, y^{2\ell})} \Sigma(b,0) \\
& - \frac{P(-y^{2b+\ell}, y^{2\ell}) P(y^{4a}, y^{2\ell})
P(y^{2a}, y^{2\ell}) P(0)^2} {P(y^{2b-2a},
y^{2\ell}) P(y^{2b+2a}, y^{2\ell}) P(y^{2b}, y^{2\ell})
P(-y^{2a+\ell}, y^{2\ell})} = 0.
\end{aligned}
\end{equation}

\noindent We now define

$$
\begin{aligned}
g(z,q) & :=z^{2} \frac{P(-q, q^{2})P(z^{4}, q^{2})}{P(z^{2}, q^{2})P(-qz^{2} ,q^{2})}\Sigma(z,1,q) - z^6 \Sigma(z^2, z, q) \\
& -  \sideset{}{'} \sum_{n=-\infty}^{\infty} \frac{(-1)^n z^{-4n} q^{n(2n+3)}}{1-q^{2n}}
\end{aligned}
$$

\noindent and

\begin{equation} \label{g}
g(a):=g(y^a, y^{\ell})=y^{2a} \frac{P(-y^{\ell}, y^{2\ell})
P(y^{4a}, y^{2\ell})}{P(y^{2a}, y^{2\ell}) P(-y^{2a+\ell}, y^{2\ell})
} \Sigma(a,0) - y^{6a} \Sigma(2a, a) - \Sigma(0, -a).
\end{equation}

\noindent The second key lemma is the following.

\begin{lem} \label{us} We have

\begin{equation} \label{part1}
\begin{aligned}
& 2g(z,q) - g(z^2, q) + 1\\
& = \frac{P(qz^2, q^2) P(-z^2, q^2) P(0)^2 P(-1, q)^2}{P(-qz^2,
q^2) P(z^2, q^2) P(-1, q^2)^2} + \frac{z^4}{2}\frac{P(q^2 z^{16},
q^4) P(-1, q^2) (q)_{\infty}^2}{P(z^8, q)}
\end{aligned}
\end{equation}

\noindent and

\begin{equation} \label{part2}
g(z,q) + g(z^{-1}q, q)=0.
\end{equation}

\end{lem}

\begin{proof}
We first require a short computation involving $\Sigma(z, \zeta,
q)$. Note that

\begin{equation} \label{sigma}
\begin{aligned}
z^4 \Sigma(z,\zeta, q) + q\zeta^4 \Sigma(zq, \zeta, q) & =
\sum_{n=-\infty}^{\infty} (-1)^n \frac{z^4 \zeta^{4n} q^{n(2n+3)}}{1 - z^2 q^{2n}} +
\sum_{n=-\infty}^{\infty} (-1)^n \frac{\zeta^{4n+4} q^{n(2n+3)}}{1-z^2 q^{2n+2}} \\
& = \sum_{n=-\infty}^{\infty} \zeta^{4n} q^{n(2n-1)} \Bigl( \frac{z^4 q^{4n} -1}{1-z^2 q^{2n}} \Bigr) \\
& = -\sum_{n=-\infty}^{\infty} (-1)^n \zeta^{4n} q^{n(2n-1)} (1 +
z^2 q^{2n})
\end{aligned}
\end{equation}

\noindent upon writing $n-1$ for $n$ in the second sum of the
first equation. Taking $\zeta=1$ yields

\begin{equation} \label{step}
z^4 \Sigma(z,1,q) + q\Sigma(zq, 1, q) = -\sum_{n=-\infty}^{\infty} (-1)^n q^{n(2n-1)} (1+ z^2 q^{2n}).
\end{equation}

\noindent Now write $g(z,q)$ in the form

$$
g(z,q)=f_{1}(z) - f_{2}(z) - f_{3}(z)
$$

\noindent where

$$
f_{1}(z):= z^2 \frac{P(-q, q^2)P(z^4,q)}{P(z^2,q^2)P(-z^2 q,q)}\Sigma(z,1,q),
$$

$$
f_{2}(z):= z^6 \Sigma(z^2, z, q),
$$

\noindent and

$$
f_{3}(z):= \sideset{}{'} \sum_{n=-\infty}^{\infty} \frac{(-1)^n
z^{-4n} q^{n(2n+3)}}{1-q^{2n}}.
$$

\noindent By (\ref{p1}) and (\ref{p2}) (replacing the base $q$ with $q^2$), and (\ref{step}),

\begin{equation} \label{f1}
f_{1}(zq) - f_{1}(z) = (z^{-2} + 1) \sum_{n=-\infty}^{\infty} (-1)^n q^{n(2n+1)}
\frac{P(z^4, q^2) P(-q, q^2)}{P(z^2,q^2) P(-z^2 q,q^2)}.
\end{equation}

\noindent A similar argument as in (\ref{sigma}) yields

\begin{equation} \label{f2}
f_{2}(zq) - f_{2}(z) = \sum_{n=-\infty}^{\infty} (-1)^n z^{4n-2}
q^{n(2n-1)} + \sum_{n=-\infty}^{\infty} (-1)^n z^{4n + 2} q^{n(2n+1)}
\end{equation}

\noindent and

\begin{equation} \label{f3}
f_{3}(zq) - f_{3}(z) = -2 + \sum_{n=-\infty}^{\infty} (-1)^n
z^{-4n} q^{n(2n-1)} + \sum_{n=-\infty}^{\infty} (-1)^n z^{-4n}
q^{n(2n+1)}.
\end{equation}

\noindent Adding (\ref{f2}) and (\ref{f3}), then subtracting
from (\ref{f1}) gives

\begin{equation} \label{constant}
g(z,q)-g(zq, q) =-2.
\end{equation}

\noindent Here we have used the identity

\begin{equation} \label{hidden1}
\begin{aligned}
& (z^{-2} + 1) \sum_{n=-\infty}^{\infty} (-1)^n q^{2n^2 + n} \frac{P(z^4, q^2) P(-q, q^2)}{P(z^2,q^2) P(-z^2 q,q^2)} \\
&= \sum_{n=-\infty}^{\infty} (-1)^n z^{4n-2}
q^{n(2n-1)} + \sum_{n=-\infty}^{\infty} (-1)^n z^{4n + 2} q^{n(2n+1)} \\
& + \sum_{n=-\infty}^{\infty} (-1)^n
z^{-4n} q^{n(2n-1)} + \sum_{n=-\infty}^{\infty} (-1)^n z^{-4n}
q^{n(2n+1)}
\end{aligned}
\end{equation}

\noindent which follows from \cite[Ex. 5.5, p.134]{Ga-Ra1}, the
triple product identity \eqref{jtp}, and a little simplification.
If we now define

$$
\begin{aligned}
f(z) & := 2g(z,q) - g(z^2, q) + 1 - \frac{P(qz^2, q^2) P(-z^2, q^2) P(0)^2 P(-1, q)^2}{P(-qz^2, q^2) P(z^2, q^2) P(-1, q^2)^2} \\
&- \frac{z^4}{2} \frac{P(q^2 z^{16}, q^4) P(-1, q^2) (q)_{\infty}^2}{P(z^8, q)},
\end{aligned}
$$

\noindent then from (\ref{p1}), (\ref{p2}), and (\ref{constant}),
one can verify that

\begin{equation} \label{functional}
f(zq) - f(z) =0.
\end{equation}

\noindent Now, it follows from a routine complex analytic argument similar to the proof of Lemma 4.2 in \cite{loveoz} (see also Lemma 2 in \cite{asd}) that $f(z)=0$. This proves (\ref{part1}).

To prove (\ref{part2}), it suffices to
show, after (\ref{constant}),

\begin{equation} \label{gees}
g(z^{-1}, q) + g(z,q)=-2.
\end{equation}

\noindent Note that

\begin{equation} \label{short}
\begin{aligned}
\Sigma(z,1, q) + z^{-6} \Sigma(z^{-1}, 1, q) & =
\sum_{n=-\infty}^{\infty} (-1)^n \frac{q^{n(2n+3)}}{1 - z^2 q^{2n}} -
z^{-4} \sum_{n=-\infty}^{\infty} (-1)^n \frac{q^{n(2n-1)}}{1-z^2 q^{2n}} \\
& = -z^{-4} \sum_{n=-\infty}^{\infty} (-1)^n q^{n(2n-1)} (1 + z^2 q^{2n})
\end{aligned}
\end{equation}

\noindent where we have written $-n$ for $n$ in the second sum in the first equation. Thus, by
(\ref{p1}), (\ref{p2}), and (\ref{short}), we have

\begin{equation} \label{f1z}
f_{1}(z) + f_{1}(z^{-1}) = -z^{-2} \sum_{n=-\infty}^{\infty} (-1)^n q^{n(2n-1)} (1+ z^2 q^{2n}) \frac{P(-q,
q^2)P(z^4,q^2)}{P(z^2,q^2)P(-z^2 q,q^2)}.
\end{equation}

\noindent Again, a similar argument as in (\ref{short}) gives

\begin{equation} \label{f2z}
f_{2}(z) + f_{2}(z^{-1}) = -z^{-2} \sum_{n=-\infty}^{\infty} (-1)^n z^{4n} q^{n(2n-1)} (1 + z^4 q^{2n})
\end{equation}

\noindent and

\begin{equation} \label{f3z}
f_{3}(z) + f_{3}(z^{-1}) = 2 - \sum_{n=-\infty}^{\infty} (-1)^n z^{4n} q^{n(2n-1)} (1+ q^{2n}).
\end{equation}

\noindent Adding (\ref{f2z}) and (\ref{f3z}), then subtracting from (\ref{f1z})
yields (\ref{gees}). Here we have used the identity

\begin{equation} \label{hidden2}
\begin{aligned}
& z^{-2} \sum_{n=-\infty}^{\infty} (-1)^n q^{n(2n-1)} (1+ z^2 q^{2n}) \frac{P(-q,
q^2)P(z^4,q^2)}{P(z^2,q^2)P(-z^2 q,q^2)} \\
&= z^{-2} \sum_{n=-\infty}^{\infty} (-1)^n z^{4n} q^{n(2n-1)} (1 + z^4 q^{2n}) +
\sum_{n=-\infty}^{\infty} (-1)^n z^{4n} q^{n(2n-1)} (1+ q^{2n})
\end{aligned}
\end{equation}

\noindent which is easily seen to be equivalent to
\eqref{hidden1}.

\end{proof}

Letting $z=y^a$ and $q=y^{\ell}$ in Lemma \ref{us}, we get

\begin{equation} \label{g1}
\begin{aligned}
& 2g(a) - g(2a) + 1 \\
&= \frac{P(y^{\ell +2a}, y^{2\ell}) P(-y^{2a}, y^{2\ell})P(0)^2
P(-1,y^{\ell})^2}{P(-y^{\ell +2a}, y^{2\ell}) P(y^{2a}, y^{2\ell})
P(-1, y^{2\ell})^2} + \frac{y^{4a}}{2} \frac{P(y^{2\ell+16a},
y^{4\ell}) P(-1,
y^{2\ell})(y^{\ell};y^{\ell})_{\infty}^2}{P(y^{8a}, y^\ell)}
\end{aligned}
\end{equation}

\noindent and

\begin{equation} \label{g2}
g(a) + g(\ell-a)=0.
\end{equation}

\noindent These two identities will be of key importance in the
proofs of Theorems \ref{main3} and \ref{main5}.


\section{Proofs of Theorems \ref{main3} and \ref{main5}}

We now compute the sums $S_{2}(3\ell -4m)$. The reason for
this choice is two-fold. First, we would like to obtain as simple
an expression as possible in the final formulation (\ref{final}).
Secondly, to prove Theorem \ref{main3}, we only need to compute $S_{2}(1)$
whereas to prove Theorem \ref{main5}, we need $S_{2}(11)$ and $S_{2}(7)$. The latter in turn yields $S_{2}(1)$ and $S_{2}(3)$ via (\ref{rels}) and
(\ref{bodd}). For $\ell=3$, we can choose $m=2$ and for $\ell=5$, $m=1$ and $m=2$ respectively. As this point, we follow the
idea of Section 6 in \cite{asd}. Namely, we write

\begin{equation} \label{n}
n=\ell r + m + b,
\end{equation}

\noindent where $-\infty < r < \infty$. The idea is to simplify
the exponent of $q$ in $S_{2}(3\ell-4m)$. Thus

$$
(3\ell-4m)n + 2n^2 = \ell^{2} r(2r+3) + 2(b+m)(b-m+\ell) + \ell(m+b) + 4b\ell r.
$$

\noindent We now substitute (\ref{n}) into (\ref{s}) and let $b$
take the values $0$, $\pm a$, and $\pm m$. Here $a$ runs through
$1$, $2$, \dotso, $\frac{\ell-1}{2}$ where the value $a \equiv \pm
m \bmod \ell$ is omitted. As in \cite{asd}, we use the notation
$\displaystyle \sideset{}{''} \sum_{a}$ to denote the sum over
these values of $a$. We thus obtain

$$
\begin{aligned}
S_{2}(3\ell-4m) & = \sideset{}{'} \sum_{n=-\infty}^{\infty} (-1)^n \frac{q^{(3\ell-4m)n + 2n^2}}{1-y^{2n}} \\
& = \sum_{b} {\sideset{}{'} \sum_{r=-\infty}^{\infty}}
(-1)^{r+m+b}  y^{m+b} q^{2(b+m)(b-m+\ell)} \frac{y^{\ell r(2r+3) + 4br}}{1-
y^{2(\ell r+m+b)}},
\end{aligned}
$$

\noindent where $b$ takes values $0$, $\pm a$, and $\pm m$ and the term corresponding
to $r=0$ and $b=-m$ is omitted. Thus

\begin{equation} \label{s(b)}
\begin{aligned}
S_{2}(3\ell-4m) & = (-1)^m y^{m} q^{2m(\ell -m)} \Sigma(m,0) + \Sigma(0, -m) + y^{6m} \Sigma(2m, m) \\
& + \sideset{}{''} \sum_{a} (-1)^{m+a} y^{m+a} q^{2(a+m)(a-m+\ell)} \Bigl
\{ \Sigma(m+a, a) + y^{-6a} \Sigma(m-a, -a) \Bigr \}.
\end{aligned}
\end{equation}

\noindent Here the first three terms arise from taking $b=0$, $-m$, and $m$ respectively.
We now can use (\ref{lem1}) to simplify this expression. By taking $b=m$ and dividing by
$y^{6a}$ in (\ref{lem1}), the sum of the two terms inside the curly brackets becomes

$$
\begin{aligned}
& y^{-4a} \frac{P(-y^{\ell}, y^{2\ell})P(y^{4a}, y^{2\ell})}{P(y^{2a}, y^{2\ell})P(-y^{2a+\ell}, y^{2\ell})} \Sigma(m,0) \\
& + y^{-6a} \frac{P(-y^{2m+\ell}, y^{2\ell}) P(y^{4a}, y^{2\ell}) P(y^{2a}, y^{2\ell}) P(0)^2}{P(y^{2m-2a}, y^{2\ell}) P(y^{2m+2a}, y^{2\ell}) P(y^{2m}, y^{2\ell}) P(-y^{2a+\ell}, y^{2\ell})}.
\end{aligned}
$$

\noindent Similarly, upon taking $a=m$ in (\ref{g}), then the sum of the second and third terms in (\ref{s(b)}) is

$$
y^{2m} \frac{P(-y^{\ell}, y^{2\ell}) P(y^{4m}, y^{2\ell})}{P(y^{2m}, y^{2\ell}) P(-y^{2m+\ell}, y^{2\ell})} \Sigma(m,0) - g(m).
$$

\noindent In total, we have

\begin{equation} \label{final}
\begin{aligned}
& S_{2}(3\ell-4m) = \\
& -g(m) \\
& + \sideset{}{''} \sum_{a} (-1)^{m+a} y^{m-5a} q^{2(a+m)(a-m+\ell)} \frac{P(-y^{2m+\ell}, y^{2\ell}) P(y^{4a}, y^{2\ell}) P(y^{2a}, y^{2\ell}) P(0)^2}{P(y^{2m-2a}, y^{2\ell}) P(y^{2m+2a}, y^{2\ell}) P(y^{2m}, y^{2\ell}) P(-y^{2a+\ell}, y^{2\ell})} \\
& + \Sigma(m,0) \Biggl\{ (-1)^m y^m q^{2m(\ell-m)} + y^{2m} \frac{P(-y^{\ell}, y^{2\ell}) P(y^{4m}, y^{2\ell})}{P(y^{2m}, y^{2\ell}) P(-y^{2m+\ell}, y^{2\ell})} \\
& + \sideset{}{''} \sum_{a} (-1)^{m+a} y^{m-3a} q^{2(a+m)(a-m+\ell)} \frac{P(-y^\ell, y^{2\ell})P(y^{4a}, y^{2\ell})}{P(y^{2a}, y^{2\ell})P(-y^{2a+\ell}, y^{2\ell})}   \Biggr\}.
\end{aligned}
\end{equation}

\noindent We can simplify some of the terms appearing in
(\ref{final}) as we are interested in certain values of $\ell$,
$m$, and $a$. To this end, we prove the following result. Let
$\{ \quad \}$ denote the coefficient of $\Sigma(m,0)$ in (\ref{final}).

\begin{prop} \label{brackets} If $\ell=3$ and $m=2$, then

$$
\{ \quad \} = -q^{9} \frac{(q^2; q^2)_{\infty} (-q^9;
q^{18})_{\infty}}{(- q; q^2)_{\infty} (q^{18}; q^{18})_{\infty}}.
$$


\noindent If $\ell=5$, $m=2$, and $a=1$, then

$$
\{ \quad \}= -q^{19} \frac{(q^2; q^2)_{\infty} (-q^{25}; q^{50})_{\infty}}{(-q; q^2)_{\infty} (q^{50}; q^{50})_{\infty}}.
$$

\noindent If $\ell=5$, $m=1$, $a=2$, then

$$
\{ \quad \}= q^{10} \frac{(q^2; q^2)_{\infty} (-q^{25}; q^{50})_{\infty}}{(-q; q^2)_{\infty} (q^{50}; q^{50})_{\infty}}.
$$

\end{prop}

\begin{proof}
These are easily deduced from Lemma \ref{lem6}.
\end{proof}

We are now in a position to prove Theorems \ref{main3} and
\ref{main5}.  We begin with Theorem \ref{main3}.

\begin{proof}
By (\ref{gf2}), (\ref{s}), and (\ref{rels}), we have

\begin{equation} \label{gen3too}
\sum_{n=0}^{\infty} \Bigl\{ N_{2} (0,3,n) -
N_{2}(1,3,n) \Bigr\} q^{n}
\frac{(q^2; q^2)_{\infty}}{(-q; q^2)_{\infty}}=2S_{2}(1) + S_{2}(7).
\end{equation}

\noindent  By (\ref{p1}), (\ref{p2}), (\ref{final}), and
Proposition \ref{brackets} we have

\begin{equation} \label{s1too}
S_{2}(1)= -g(2) - y^3 \frac{(q^2; q^2)_{\infty} (-q^9;
q^{18})_{\infty}}{(-q; q^2)_{\infty} (q^{18}; q^{18})_{\infty}}
\Sigma(2,0).
\end{equation}

\noindent Taking $b=1$ in (\ref{bodd}) yields
\begin{equation}
S_{2}(1) - S_{2}(7) = \frac{(q^2; q^2)_{\infty}}{(-q; q^2)_{\infty}} -1.
\end{equation}

We need to prove that

$$
\begin{aligned}
 & -3g(2) - 3y^2 \frac{(q^2; q^2)_{\infty} (-q^9; q^{18})_{\infty}}{(-q; q^2)_{\infty} (q^{18}; q^{18})_{\infty}} \Sigma(2,0) - \frac{(q^2; q^2)_{\infty}}{(-q; q^2)_{\infty}} +1\\
 &= \Bigl \{r_{01}(0) q^0 + r_{01}(1)q + r_{01}(2)q^2  \Bigr \} \frac{(q^2; q^2)_{\infty}}{(-q; q^2)_{\infty}}.
 \end{aligned}
$$

\noindent We now multiply the right hand side of the above
expression using Lemma \ref{lem6} and the $r_{01}(d)$ from Theorem
\ref{main3} (recall that $r_{01}(d)$ is just $R_{01}(d)$ with $q$
replaced by $q^{3}$).  We then equate coefficients of powers of
$q$ and verify the resulting identities.  The only power of $q$
for which the resulting equation does not follow easily upon
cancelling factors in infinite products is the constant term. We
obtain

$$
\begin{aligned}
-3g(2) + 1 & = \frac{(q^{18}; q^{18})_{\infty}^4 (-q^9; q^9)_{\infty}^4 (q^3; q^6)_{\infty} (q^3, -q^6, -q^9, -q^{12}, q^{15}, q^{18}; q^{18})_{\infty}}{(q^{12}; q^{12})_{\infty} (q^6, q^{30}, q^{36}; q^{36})_{\infty}^2} \\
& - y \frac{(q^9; q^9)_{\infty} (-q^{18}; q^{18})_{\infty} (q^9, q^{27}, q^{36}; q^{36})_{\infty}}{(q^3, q^{15}; q^{18})_{\infty} (q^{12}, q^{24}; q^{36})_{\infty}}.
\end{aligned}
$$

\noindent The first term (resp. second term) above is easily seen to be
identical to the first term (resp. second term) in \eqref{g1} with
$a=1$.  Applying \eqref{g2}, this then establishes the above
identity and completes the proof of Theorem \ref{main3}.
\end{proof}

We now turn to Theorem \ref{main5}.

\begin{proof}
We begin with the rank differences $R_{12}(d)$.  By (\ref{gf2}),
(\ref{s}), and (\ref{rels}), we have

\begin{equation} \label{gen3}
\sum_{n=0}^{\infty} \Bigl\{ N_{2}(1,5,n) -
N_{2}(2,5,n) \Bigr\} q^{n} \frac{(q^2; q^2)_{\infty}}{(-q; q^2)_{\infty}} = 2S_{2}(3) - S_{2}(1).
\end{equation}

\noindent and by (\ref{p1}), (\ref{p2}), (\ref{final}), and Proposition \ref{brackets},

\begin{equation} \label{s1}
S_{2}(1)=  -g(1) + q \frac{P(0)^2 P(-y^7, y^{10})}{P(y^2, y^{10}) P(-y^9, y^{10})} + y^2 \Sigma(1,0) \frac{(q^2; q^2)_{\infty} (-q^{25}; q^{50})_{\infty}}{(-q; q^2)_{\infty} (q^{50}; q^{50})_{\infty}} + \frac{(q^2; q^2)_{\infty}}{(-q; q^2)_{\infty}} - 1
\end{equation}

\noindent and

\begin{equation} \label{s3}
S_{2}(3) =  g(2) + yq^4 \frac{P(0)^2 P(-y^9, y^{10})}{P(y^4, y^{10}) P(-y^7, y^{10})} + y^3 q^4 \Sigma(2,0) \frac{(q^2; q^2)_{\infty} (-q^{25}; q^{50})_{\infty}}{(-q; q^2)_{\infty} (q^{50}; q^{50})_{\infty}}.
\end{equation}

\noindent By (\ref{gen3}), (\ref{s1}), and (\ref{s3}), we need to prove

$$
\begin{aligned}
& 2g(2) + 2yq^4 \frac{P(0)^2 P(-y^9, y^{10})}{P(y^4, y^{10}) P(-y^7, y^{10})} + 2y^3 q^4 \Sigma(2,0) \frac{(q^2; q^2)_{\infty} (-q^{25}; q^{50})_{\infty}}{(-q; q^2)_{\infty} (q^{50}; q^{50})_{\infty}}\\
& g(1) - q \frac{P(0)^2 P(-y^2, y^{10})}{P(y^2, y^{10}) P(-y^9, y^{10})} - y^2 \Sigma(1,0) \frac{(q^2; q^2)_{\infty} (-q^{25}; q^{50})_{\infty}}{(-q; q^2)_{\infty} (q^{50}; q^{50})_{\infty}} - \frac{(q^2; q^2)_{\infty}}{(-q; q^2)_{\infty}} + 1 \\
& = \Bigl \{ r_{12}(0) q^0 + r_{12}(1)q + r_{12}(2)q^2 +
r_{12}(3)q^3 + r_{12}(4)q^4 \Bigr \}
\frac{(q^2; q^2)_{\infty}}{(-q; q^2)_{\infty}}.
\end{aligned}
$$

\noindent We now multiply the right hand side of the above
expression using Lemma \ref{lem6} and the $R_{12}(d)$ from Theorem
\ref{main5}, equating coefficients of powers of $q$.  The
coefficients of $q^{0}$, $q^{1}$, $q^{2}$, $q^{3}$, $q^{4}$ give
us, respectively,

\begin{equation} \label{check0}
\begin{aligned}
& 2g(2) + g(1) + 1\\
& = \frac{(q^5, q^{45}; q^{50})_{\infty}^2 (q^{30}, q^{40}, q^{60}, q^{70}; q^{100})_{\infty} (q^{50}; q^{100})_{\infty}^3 (q^{100}; q^{100})_{\infty}^2 (-q^{10}, q^{15}, -q^{25}, q^{35}, -q^{40}, q^{50}; q^{50})_{\infty}}{(q^5; q^5)_{\infty}} \\
& - y \frac{(q^{15}, q^{35}, q^{50}; q^{50})_{\infty}^2 (q^5, -q^{20}, -q^{25}, -q^{30}, q^{45}, q^{50}; q^{50})_{\infty}}{(q^5; q^{10})_{\infty} (q^{30}, q^{40}, q^{60}, q^{70}, q^{100}; q^{100})_{\infty}} \\
& - y^2 \frac{(q^{10}, q^{90}; q^{100})_{\infty} (q^{25}; q^{25})_{\infty} (-q^{50}; q^{50})_{\infty} (q^{25}, q^{75}, q^{100}; q^{100})_{\infty}}{(q^5, q^{20}; q^{25})_{\infty}},
\end{aligned}
\end{equation}

\begin{equation} \label{check1}
\begin{aligned}
& \frac{(q^{50}; q^{50})_{\infty}^2 (-q^{15}, -q^{35}; q^{50})_{\infty}}{(q^{10}, q^{40}; q^{50})_{\infty} (-q^5, -q^{45}; q^{50})_{\infty}} \\
& = \frac{(q^5, q^{45}; q^{50})_{\infty}^2 (q^{30}, q^{40}, q^{60}, q^{70}; q^{100})_{\infty} (q^{50}; q^{100})_{\infty}^3 (q^{100}; q^{100})_{\infty}^2 (q^{5}, -q^{20}, -q^{25}, -q^{30}, q^{45}, q^{50}; q^{50})_{\infty}}{(q^5; q^5)_{\infty}} \\
& + y \frac{(-q^{25}, q^{50}; q^{50})_{\infty} (q^{25}, q^{75}, q^{100}; q^{100})_{\infty}}{(q^{20}, q^{30}; q^{50})_{\infty}},
\end{aligned}
\end{equation}

\begin{equation} \label{check2}
\begin{aligned}
& \frac{(q^{10}, q^{90}; q^{100})_{\infty} (q^{25}; q^{25})_{\infty} (-q^{50}; q^{50})_{\infty} (-q^{10}, q^{15}, -q^{25}, q^{35}, -q^{40}, q^{50}; q^{50})_{\infty}}{(q^5, q^{20}; q^{25})_{\infty}}\\
& = \frac{(q^{15}, q^{35}, q^{50}; q^{50})_{\infty}^2 (q^{25}, q^{75}, q^{100}; q^{100})_{\infty}}{(q^5; q^{10})_{\infty} (q^{30}, q^{40}, q^{60}, q^{70}, q^{100}; q^{100})_{\infty}},
\end{aligned}
\end{equation}

\begin{equation} \label{check3}
 \begin{aligned}
& \frac{(-q^{25}, q^{50}; q^{50})_{\infty} (-q^{10}, q^{15}, -q^{25}, q^{35}, -q^{40}, q^{50}; q^{50})_{\infty}}{(q^{20}, q^{30}; q^{50})_{\infty}} \\
 & = y\frac{(q^{10}, q^{90}; q^{100})_{\infty} (q^{25}; q^{25})_{\infty} (-q^{50}; q^{50})_{\infty} (q^5, -q^{20}, -q^{25}, -q^{30}, q^{45}, q^{50}; q^{50})_{\infty}}{(q^5, q^{20}; q^{25})_{\infty}} \\
 & +  \frac{(q^5, q^{45}; q^{50})_{\infty}^2 (q^{30}, q^{40}, q^{60}, q^{70}; q^{100})_{\infty} (q^{50}; q^{100})_{\infty}^3 (q^{100}; q^{100})_{\infty}^2 (q^{25}, q^{75}, q^{100}; q^{100})_{\infty}}{(q^5; q^5)_{\infty}},
\end{aligned}
\end{equation}

\begin{equation} \label{check4}
\begin{aligned}
2y \frac{(q^{50}; q^{50})_{\infty}^2 (-q^5, -q^{45}; q^{50})_{\infty}}{(q^{20}, q^{30}; q^{50})_{\infty} (-q^{15}, -q^{35}; q^{50})_{\infty}} & = \frac{(q^{15}, q^{35}, q^{50}; q^{50})_{\infty}^2 (-q^{10}, q^{15}, -q^{25}, q^{35}, -q^{40}, q^{50}; q^{50})_{\infty}}{(q^5; q^{10})_{\infty} (q^{30}, q^{40}, q^{60}, q^{70}, q^{100}; q^{100})_{\infty}} \\
& -  \frac{(-q^{25}, q^{50}; q^{50})_{\infty} (q^5, -q^{20}, -q^{25}, -q^{30}, q^{45}, q^{50}; q^{50})_{\infty}  }{(q^{20}, q^{30}; q^{50})_{\infty}}.
\end{aligned}
\end{equation}

Equation \eqref{check2} is immediate after some simplification.
The other identities follow from routine (though tedious)
reduction and application of one of Lemmas \ref{Hick1} -
\ref{Hick2.5}. Specifically, upon clearing denominators in
\eqref{check1} and simplifying, we have
\begin{equation}
\begin{aligned}
(-q^{15},q^{20},q^{30},-q^{35}; q^{50})_{\infty} &=
(q^5,q^{10},-q^{15},-q^{20},-q^{25},-q^{25},-q^{30},-q^{35},q^{40},q^{45};q^{50})_{\infty}
\\ &+ y(-q^5,q^{10},q^{40},-q^{45};q^{50})_{\infty}.
\end{aligned}
\end{equation}
Now replacing $q$ by $-q$, this may be verified using the case
$(x,z,q) = (-q^5,q^{10},q^{25})$ of Lemma \ref{Hick1}.  After
clearing denominators and simplifying, \eqref{check3} may be
reduced to

\begin{equation} \label{yikes}
\begin{aligned}
(-q^{25};q^{50})_{\infty}^2(q^{20},q^{80};q^{100})_{\infty}(q^{15},q^{35};q^{50})_{\infty}
&=
y(q^{10},q^{90};q^{100})_{\infty}(q^{10},-q^{20},-q^{30},q^{40};q^{50})
\\ &+
(q^{40},q^{60};q^{100})_{\infty}(q^5,-q^{15},-q^{35},q^{45};q^{50})_{\infty}.
\end{aligned}
\end{equation}

\noindent Factoring out $(q^5,-q^{20},-q^{30},q^{45};q^{50})_{\infty}$
from the right hand side, replacing $q$ by $-q$ and applying the
case $(x,z,q) = (q^5,-q^{10},q^{25})$ of Lemma \ref{Hick1}
verifies (\ref{yikes}).  For \eqref{check4}, we clear denominators
and simplify to get

\begin{equation} \label{wow}
\begin{aligned}
& 2y(q^{10},q^{25},q^{40};q^{50})_{\infty}(q^{100};q^{100})_{\infty} \\
&= (q^{10},-q^{10},-q^{10},q^{15},q^{15},-q^{15},-q^{25},q^{35},q^{35},-q^{35},q^{40},-q^{40},-q^{40},q^{50};q^{50})_{\infty}
\\
& -
(q^{5},q^{5},-q^{15},q^{20},-q^{20},-q^{20},-q^{25},q^{30},-q^{30},-q^{30},-q^{35},q^{45},q^{45},q^{50};q^{50})_{\infty}.
\end{aligned}
\end{equation}

\noindent Factoring out $(-q^{15},-q^{25},-q^{35},q^{50};q^{50})_{\infty}$
from both terms on the right hand side, replacing $q$ by $-q$ and
writing the right hand side in base $q^{25}$ yields an expression
to which the case $(x,z,q) = (-q^5,-q^{10},q^{25})$ of Lemma
\ref{Hick2} may applied, confirming (\ref{wow}).


As for \eqref{check0}, taking $a=2$ in \eqref{g1} and applying
\eqref{g2} gives
\begin{equation}
\begin{aligned}
2g(2)+g(1)+1 &=
\frac{(q^5,q^{45};q^{50})_{\infty}(-q^{20},-q^{30};q^{50})_{\infty}
(-q^{25};q^{50})_{\infty}^4(q^{50};q^{50})_{\infty}^2}
{(-q^5,-q^{45};q^{50})_{\infty}(q^{20},q^{30};q^{50})_{\infty}}
\\ &-
y^2\frac{(q^{10},q^{90};q^{100})_{\infty}(q^{25};q^{25})_{\infty}^2(-q^{50};q^{50})_{\infty}^2}
{(q^5,q^{20};q^{25})_{\infty}}.
\end{aligned}
\end{equation}
The final term above is identical to the final term in
\eqref{check0}.  After some simplification, the fact the the first
term above is equal to the first two terms in \eqref{check0} is
equivalent to the identity

\begin{equation} \label{goodness}
\begin{aligned}
(q^{30},q^{70};q^{100})_{\infty}(-q^{10},-q^{15},-q^{35},-q^{40};q^{50})_{\infty}
&-y(-q^5,-q^{45};q^{50})_{\infty} \\ &=
(q^5,-q^{15},-q^{25},-q^{25},-q^{30},-q^{35},q^{45};q^{50})_{\infty}.
\end{aligned}
\end{equation}

\noindent Equation (\ref{goodness}) is seen to be true after multiplying both sides by
$(q^{10},q^{40};q^{50})_{\infty}$, replacing $q$ by $-q$, and
applying the case $(x,z,q) = (-q^5,q^{10},q^{25})$ of Lemma
\ref{Hick1}.

We now turn to the rank differences $R_{02}(d)$, proceeding as
above. Again by (\ref{gf2}), (\ref{s}), and (\ref{rels}), we have

\begin{equation} \label{gen4}
\sum_{n=0}^{\infty} \Bigl\{ N_{2}(0,5,n) -
N_{2}(2,5,n) \Bigr\} q^{n} \frac{(q^2; q^2)_{\infty}}{(-q; q^2)_{\infty}} = 2S_{2}(1) + S_{2}(3) - \frac{(q^2; q^2)_{\infty}}{(-q; q^2)_{\infty}} + 1.
\end{equation}

\noindent By (\ref{gen4}),
(\ref{s1}), and (\ref{s3}), it suffices to prove

$$
\begin{aligned}
&   -2g(1) + 2q \frac{P(0)^2 P(-y^7, y^{10})}{P(y^2, y^{10}) P(-y^9, y^{10})} + 2y^2 \Sigma(1,0) \frac{(-q^{25}; q^{50})_{\infty}}{(-q; q^2)_{\infty} (q^{50}; q^{50})_{\infty}} \\
&  + g(2) + yq^4 \frac{P(0)^2 P(-y^9, y^{10})}{P(y^4, y^{10}) P(-y^7, y^{10})} + y^3q^4 \Sigma(2,0) \frac{(-q^{25}; q^{50})_{\infty}}{(-q; q^2)_{\infty} (q^{50}; q^{50})_{\infty}} + \frac{(q^2; q^2)_{\infty}}{(-q; q^2)_{\infty}} - 1. \\
& = \Bigl \{ r_{02}(0) q^0 + r_{02}(1)q + r_{02}(2)q^2 +
r_{02}(3)q^3 + r_{02}(4)q^4 \Bigr \} \frac{(q^2; q^2)_{\infty}}{(-q; q^2)_{\infty}}.
 \end{aligned}
$$

\noindent Again, equating coefficients of powers of $q$ yields the
following identities.

\begin{equation} \label{check5}
\begin{aligned}
& 2g(1) - g(2) +1\\
&= \frac{(q^5, q^{45}; q^{50})_{\infty}^2 (q^{50}; q^{50})_{\infty}^3 (q^{30}, q^{40}, q^{60}, q^{70}; q^{100})_{\infty} (-q^{10}, q^{15}, -q^{25}, q^{35}, -q^{40}, q^{50}; q^{50})_{\infty}}{(q^5; q^5)_{\infty} (q^{100}; q^{100})_{\infty}} \\
& + y \frac{(q^{15}, q^{35}, q^{50}; q^{50})_{\infty}^2 (q^5, -q^{20}, -q^{25}, -q^{30}, q^{45}, q^{50}; q^{50})_{\infty}}{(q^5; q^{10})_{\infty} (q^{30}, q^{40}, q^{60}, q^{70}, q^{100}; q^{100})_{\infty}} \\
& + y \frac{(q^{30}, q^{70}; q^{100})_{\infty} (q^{25}; q^{25})_{\infty} (-q^{50}; q^{50})_{\infty} (q^{25}, q^{75}, q^{100}; q^{100})_{\infty}}{(q^{10}, q^{15}; q^{25})_{\infty}},
\end{aligned}
\end{equation}

\begin{equation} \label{check6}
\begin{aligned}
& 2\frac{(q^{50}; q^{50})_{\infty}^2 (-q^{15}, -q^{35}; q^{50})_{\infty}}{(q^{10}, q^{40}; q^{50})_{\infty} (-q^5, -q^{45}; q^{50})_{\infty}} \\
& = \frac{(-q^{25}; q^{50})_{\infty} (q^{50}; q^{50})_{\infty} (-q^{10}, q^{15}, -q^{25}, q^{35}, -q^{40}, q^{50}; q^{50})_{\infty}}{(q^{10}, q^{40}; q^{50})_{\infty}} \\
& + \frac{(q^5, q^{45}; q^{50})_{\infty}^2 (q^{50}; q^{50})_{\infty}^3 (q^{30}, q^{40}, q^{60}, q^{70}; q^{100})_{\infty} (q^5, -q^{20}, -q^{25}, -q^{30}, q^{45}, q^{50}; q^{50})_{\infty}}{(q^5; q^5)_{\infty} (q^{100}; q^{100})_{\infty}},
\end{aligned}
\end{equation}

\begin{equation} \label{check7}
\begin{aligned}
& \frac{(q^{30}, q^{70}; q^{100})_{\infty} (q^{25}; q^{25})_{\infty} (-q^{50}; q^{50})_{\infty} (-q^{10}, q^{15}, -q^{25}, q^{35}, -q^{40}, q^{50}; q^{50})_{\infty}}{(q^{10}, q^{15}; q^{25})_{\infty}} \\
& = \frac{(-q^{25}; q^{50})_{\infty} (q^{50}; q^{50})_{\infty} (q^5, -q^{20}, -q^{25}, -q^{30}, q^{45}, q^{50}; q^{50})_{\infty}}{(q^{10}, q^{40}; q^{50})_{\infty}} \\
& + y  \frac{(q^{15}, q^{35}, q^{50}; q^{50})_{\infty}^2 (q^{25},
q^{75}, q^{100}; q^{100})_{\infty}}{(q^5; q^{10})_{\infty}
(q^{30}, q^{40}, q^{60}, q^{70}, q^{100}; q^{100})_{\infty}},
\end{aligned}
\end{equation}

\begin{equation} \label{check8}
\begin{aligned}
&\frac{(q^{30}, q^{70}; q^{100})_{\infty} (q^{25}; q^{25})_{\infty} (-q^{50}; q^{50})_{\infty} (q^5, -q^{20}, -q^{25}, -q^{30}, q^{45}, q^{50}; q^{50})_{\infty}}{(q^{10}, q^{15}; q^{25})_{\infty}} \\
& = \frac{(q^5, q^{45}; q^{50})_{\infty}^2 (q^{50}; q^{50})_{\infty}^3 (q^{30}, q^{40}, q^{60}, q^{70}; q^{100})_{\infty} (q^{25}, q^{75}, q^{100}; q^{100})_{\infty} }{(q^5; q^5)_{\infty} (q^{100}; q^{100})_{\infty}},
\end{aligned}
\end{equation}

\begin{equation} \label{check9}
\begin{aligned}
& y\frac{(q^{50}; q^{50})_{\infty}^2 (-q^5, -q^{45}; q^{50})_{\infty}}{(q^{20}, q^{30}; q^{50})_{\infty} (-q^{15}, -q^{35}; q^{50})_{\infty}} \\
& = \frac{(q^{15}, q^{35}, q^{50}; q^{50})_{\infty}^2 (-q^{10}, q^{15}, -q^{25}, q^{35}, -q^{40}, q^{50}; q^{50})_{\infty}}{(q^5; q^{10})_{\infty} (q^{30}, q^{40}, q^{60}, q^{70}, q^{100}; q^{100})_{\infty}} \\
& -  \frac{(-q^{25}; q^{50})_{\infty} (q^{50}; q^{50})_{\infty} (q^{25}, q^{75}, q^{100}; q^{100})_{\infty}}{(q^{10}, q^{40}; q^{50})_{\infty}}.
\end{aligned}
\end{equation}

These follow in the same way as equations \eqref{check0} -
\eqref{check4}.  The arduous details are left to the
interested reader. The point is to simplify and reduce in order to arrive at
an expression that can be verified using an appropriate instance
of one of the Lemmas \ref{Hick1} - \ref{Hick2.5}.



\end{proof}

\section{Concluding Remarks}
With the present paper and previous work on rank differences for
overpartitions \cite{loveoz}, we have seen the effectiveness of
the approach developed by Atkin and Swinnerton-Dyer \cite{asd} for
proving formulas for rank differences in arithmetic progressions
in terms of modular forms and generalized Lambert series. We
should stress that two major difficulties in this method are the
requirement that all of the formulas be ascertained beforehand and
the apparent need for a new set of key $q$-series identities for
each application. Nevertheless, the ideas should in principle be
reliable in other instances where there is a two-variable
generating function like \eqref{gf1.5}.  For example, one might
consider the ranks arising in Andrews' study of Durfee symbols
\cite{An2} or the generalized ranks of Garvan \cite{Ga1}. Finally,
as evidenced by work of Atkin and Hussain \cite{At-Hu1} on the partition rank,
the formulas for rank differences quickly become more complicated
as $\ell$ grows.  It would be interesting to try to extend the
method used for $\ell = 3$ and $5$ here and in \cite{loveoz} to
the case $\ell = 7$.

\section*{Acknowlegements}
The second author would like to thank the Institut des Hautes {\'E}tudes
Scientifiques for their hospitality and support during the preparation
of this paper.

\end{document}